\documentclass[preprint,12pt]{elsarticle}

\usepackage{amsmath,amssymb,amsthm,mathtools,mathrsfs}
\usepackage{enumitem,microtype,tikz}
\usepackage[colorlinks=true,linkcolor=blue,citecolor=blue,urlcolor=blue]{hyperref}
\usetikzlibrary{arrows.meta,positioning}
\journal{Topology and its Applications}
\myfooter[L]{Preprint submitted to Topology and its Applications}
\biboptions{numbers,sort&compress}
\allowdisplaybreaks[2]

\newtheorem{theorem}{Theorem}[section]
\newtheorem{lemma}[theorem]{Lemma}
\newtheorem{proposition}[theorem]{Proposition}
\newtheorem{corollary}[theorem]{Corollary}
\newtheorem{definition}[theorem]{Definition}
\newtheorem{example}[theorem]{Example}

\begin{document}
\begin{frontmatter}
\title{Strongly meet continuous posets}
\author[inst1]{Yifeng Xie}
\author[inst1]{Bin Zhao\corref{cor1}}
\ead{zhaobin@snnu.edu.cn}
\affiliation[inst1]{organization={School of Mathematics and Statistics, Shaanxi Normal University},
  city={Xi'an}, postcode={710119}, country={China}}
\cortext[cor1]{Corresponding author.}

\begin{abstract}
The main purpose of this paper is to introduce and investigate strongly meet continuous posets. It is shown that a poset is strongly meet continuous if and only if the Scott closure of the set of elements way below each element coincides with its principal ideal. Moreover, the Scott topology of every strongly meet continuous poset is shown to have a basis consisting of Scott open filters, and a countable strongly meet continuous dcpo that is not continuous is constructed. For posets with the one step closure property, strongly meet continuity is equivalent to continuity. Further, a poset is continuous if and only if it is strongly meet continuous and its Scott space is core compact. As a consequence, every core compact strongly meet continuous dcpo has a sober Scott space.
\end{abstract}

\end{frontmatter}

\section{Introduction}\label{sec:introduction}

The Scott topology and the way below relation play a fundamental role in domain
theory\cite{ref18,ref11}. In particular, continuous posets admit approximation by directed sets of elements way below a given point, and this order theoretic property is
closely connected with the structure of their Scott topologies. It is natural to study conditions under which a poset can be recovered from its approximating elements and to seek topological characterizations of continuity.

Meet continuity provides a weakening of continuity formulated in terms of
Scott closure and directed suprema\cite{ref12}.
Another approach to continuity is through the one step closure property, which
expresses Scott closure in terms of directed suprema. The relationships among
one step closure, weak one step closure and meet continuity have been studied
in \cite{ref4,ref5}. In particular, A poset have one step closure can not be
continuous. And a poset is continuous if and only if it is both
 quasicontinuous and meet continuous \cite{ref19}. These results lead to the
question of whether exist the posets that strengthening of meet continuous, and still weaker than continuous. It can be combined with suitable order theoretic or topological
conditions to characterize continuous posets.

In this paper, we consider strongly meet continuous poset. It follows that every continuous poset is strongly meet continuous and every strongly meet continuous poset is meet continuous. We construct a countable
strongly meet continuous dcpo which is not continuous. We then examine which additional conditions can make strongly meet continuous poset become continuous poset. For strongly
meet continuous posets, we prove that one step closure and weak one step closure
are each equivalent to continuity. On the topological side, we show that
Scott open filters form a basis for the Scott topology of every strongly meet
continuous poset. This yields a characterization of continuous posets by core
compactness of their Scott spaces and the existence of such a filter basis;
in particular, a strongly meet continuous poset is continuous if and only if
its Scott space is core compact.

The paper is organized as follows. In Section~\ref{sec:preliminaries} we recall
some basic notions and results. In Section~\ref{sec:strongly-meet} we introduce
strongly meet continuous poset, establish its basic properties and give the countable
dcpo example. Section~\ref{sec:continuous-case} is devoted to characterizations
of continuity within the class of strongly meet continuous posets.
Finally, Section~\ref{sec:conclusion} summarizes the main results and points out a direction for future research.

\section{Preliminaries}\label{sec:preliminaries}

We recall some notions and results from order theory and topology used in this paper. Further details can be found in \cite{ref11,ref14,ref18}.

Let $(P,\leq)$ be a poset. For $A\subseteq P$, put
$\downarrow A=\{x\in P:\exists\ a\in A,\ x\leq a\}$ and
$\uparrow A=\{x\in P:\exists\ a\in A,\ a\leq x\}$.
A nonempty subset $D\subseteq P$ is \emph{directed} if, for every finite subset $F\subseteq D$, there exists $d\in D$ such that $f\leq d$ for all $f\in F$. Let
$\mathcal D^\vee(P)=\{D\subseteq P:D\ne\varnothing,\ D\text{ is directed and }\bigvee D\text{ exists}\}$.
If every nonempty directed subset of $P$ has a supremum, then $P$ is a \emph{directed complete partial order} (dcpo).
A subset $U\subseteq P$ is \emph{Scott open} if
$U=\uparrow U$, and for every $D\in\mathcal D^\vee(P)$,
$\bigvee D\in U$ implies $D\cap U\ne\varnothing$.
The topology of Scott open sets is denoted by $\sigma(P)$, and the corresponding space by $\Sigma P$. We write $\Gamma(P)$ for the family of Scott closed subsets and $\operatorname{cl}_{\sigma}$ and $\operatorname{int}_{\sigma}$ for the Scott closure and interior operators, respectively. In particular,
$\operatorname{cl}_{\sigma}\{x\}=\downarrow x$ and
$\operatorname{cl}_{\sigma}(\downarrow A)=\operatorname{cl}_{\sigma}A$.
The \emph{way below relation} is denoted by
$a\ll x$ if and only if $\forall\ D\in\mathcal D^\vee(P), x\leq\bigvee D\ ,$ then $\exists\ d\in D,\ a\leq d.$
Put
$\Downarrow x=\{a\in P:a\ll x\}$ and
$\Uparrow a=\{x\in P:a\ll x\}$.
A poset $P$ is \emph{continuous} if for every $x\in P$, the set $\Downarrow x$ is directed and $\bigvee\Downarrow x=x$. A continuous dcpo is called a \emph{domain}.

\begin{lemma}[\textup{\cite{ref11}}]\label{lem:directed-approx}
	Let $P$ be a poset. If, for every $x\in P$, there exists a directed set $E_x\subseteq\ \Downarrow x$ such that $\bigvee E_x=x$, then $P$ is continuous.
\end{lemma}

\begin{lemma}[\textup{\cite{ref11}}]\label{lem:interpolation}
	If $P$ is continuous, then the way below relation has the interpolation property.
	Moreover, $\{\Uparrow a:a\in P\}$ is a base for $\sigma(P)$.
\end{lemma}

A poset $P$ is \emph{meet continuous} if
$\forall\ x\in P,\ \forall\ D\in\mathcal D^\vee(P),x\leq\bigvee D,$ then $
x\in\operatorname{cl}_{\sigma}(\,\downarrow x\ \cap\downarrow D).$
For a topological space $X$, let $\mathcal O(X)$ be its lattice of open sets. For $V,U\in\mathcal O(X)$, define
$V\ll U$ if and only if every open cover of $U$ has a finite subfamily covering $V$.
The space $X$ is \emph{core compact} if
$\forall\ x\in U\in\mathcal O(X), \exists\ V\in\mathcal O(X), x\in V\ll U.$
A nonempty subset $A$ of a topological space $X$ is \emph{irreducible} if
$\forall\, U,V\in\mathcal O(X), A\cap U\neq\varnothing$ and $\ A\cap V\neq\varnothing,$ then $
A\cap U\cap V\neq\varnothing.$
A $T_0$ space is \emph{sober} if every nonempty irreducible closed subset is the closure of a unique point.
A \emph{Scott open filter} is a nonempty, lower directed, Scott open subset of $P$.

\section{Strongly meet continuous posets}\label{sec:strongly-meet}
In this section, we introduce strongly meet continuous posets and investigate its basic properties. We examine its relationship with continuity and meet continuity, and construct a countable strongly meet continuous dcpo which is not continuous. 
\begin{definition}\label{def:smc}\normalfont
A poset $P$ is called \emph{strongly meet continuous} if, for every $x\in P$ and every family $\mathscr D\subseteq\mathcal D^{\vee}(P)$ satisfying $x\leq \vee D$ for every $D\in \mathscr D$, we have 
 $x\in\operatorname{cl}_{\sigma}\!\left(\mathord{\downarrow} x\ \cap\ \bigcap_{D\in\mathscr D}\mathord{\downarrow} D\right).$
\end{definition}

\begin{theorem}\label{thm:smc-equivalences}
Let $P$ be a poset, the following conditions are equivalent:
\begin{enumerate}[label=(\roman*)]
 \item $P$ is strongly meet continuous;
 \item $x\in\operatorname{cl}_{\sigma}(\mathord{\Downarrow} x)$ for every $x\in P$;
 \item $\operatorname{cl}_{\sigma}(\mathord{\Downarrow} x)=\mathord{\downarrow} x$ for every $x\in P$;
 \item $x\in U$ implies $U\cap\mathord{\Downarrow} x\ne\varnothing$ for every $x\in P$ and $U\in\sigma(P)$.
\end{enumerate}
\end{theorem}
\begin{proof}
For each $x\in P$, put
$\mathscr D_x=\{D\in\mathcal D^{\vee}(P):x\leq\bigvee D\}$.
We first claim that
\[
 \mathord{\Downarrow} x=\mathord{\downarrow} x\ \cap\ \bigcap_{D\in\mathscr D_x}\mathord{\downarrow} D. \tag{$\ast$}
\]
Indeed, if $a\ll x$, then $a\leq x$, and for every $D\in\mathscr D_x$ there exists
$d\in D$ with $a\leq d$, so $a\in\mathord{\downarrow} x\cap\bigcap_{D\in\mathscr D_x}\mathord{\downarrow} D$.
Conversely, if $a\in \mathord{\downarrow} x\ \cap\ \bigcap_{D\in\mathscr D_x}\mathord{\downarrow} D,$ and $D\in\mathcal D^{\vee}(P)$ satisfies
$x\leq\bigvee D$, then $D\in\mathscr D_x$, and hence $a\in\mathord{\downarrow} D$; this means that
$a\leq d$ for some $d\in D$, which proves $a\ll x$.

(i)$\Rightarrow$(ii). Applying Definition~\ref{def:smc} to the family $\mathscr D_x$ and
using $(\ast)$, we obtain
$x\in\operatorname{cl}_{\sigma}(\mathord{\downarrow} x\cap\bigcap_{D\in\mathscr D_x}\mathord{\downarrow} D)=\operatorname{cl}_{\sigma}(\mathord{\Downarrow} x)$.

(ii)$\Rightarrow$(i). Let $\mathscr D\subseteq\mathcal D^{\vee}(P)$ that $x\leq\bigvee D$
for every $D\in\mathscr D$. Then $\mathscr D\subseteq\mathscr D_x$, so by $(\ast)$,
$\mathord{\Downarrow} x\subseteq\mathord{\downarrow} x\ \cap\ \bigcap_{D\in\mathscr D}\mathord{\downarrow} D.$
Taking Scott closures on both sides and applying (ii), we conclude that
\[
 x\in\operatorname{cl}_{\sigma}(\mathord{\Downarrow} x)\subseteq
 \operatorname{cl}_{\sigma}\!\left(\mathord{\downarrow} x\ \cap\ \bigcap_{D\in\mathscr D}\mathord{\downarrow} D\right).
\]

(ii)$\Leftrightarrow$(iii). Since $\mathord{\Downarrow} x\subseteq\mathord{\downarrow} x$ and $\mathord{\downarrow} x$ is Scott
closed, we always have $\operatorname{cl}_{\sigma}(\mathord{\Downarrow} x)\subseteq\mathord{\downarrow} x$. Hence the equality
$\operatorname{cl}_{\sigma}(\mathord{\Downarrow} x)=\mathord{\downarrow} x$ holds if and only if $x\in\operatorname{cl}_{\sigma}(\mathord{\Downarrow} x)$.

(ii)$\Leftrightarrow$(iv). This follows from the standard neighborhood characterization
of closure.
\end{proof}

\begin{proposition}\label{prop:continuous-smc-mc}
	Every continuous poset is strongly meet continuous, and every strongly meet continuous
	poset is meet continuous.
\end{proposition}
\begin{proof}
	Suppose that $P$ is continuous. For each $x\in P$, the set $\mathord{\Downarrow} x$ is directed and
	$\bigvee\mathord{\Downarrow} x=x$. So $P$ is strongly meet
	continuous by Theorem~\ref{thm:smc-equivalences}. And if $P$ is strongly meet
	continuous, then taking $\mathscr D=\{D\}$ in Definition~\ref{def:smc} shows precisely
	that $P$ is meet continuous.
\end{proof}

In the following we construct a countable strongly meet continuous dcpo that is not continuous. This shows that the converse
of the first implication in Proposition~\ref{prop:continuous-smc-mc} fails.

Let
\[
{P_{\mathrm w}}=\{\varepsilon\}\cup\bigcup_{k\geq 1}\mathbb N_{+}^k,
\]
where $\varepsilon$ is the empty sequence of length zero. We denote the length of a
sequence $s$ by $|s|$, and the first and last entries of a nonempty sequence $s$ by
$\operatorname{f}(s)$ and $v(s)$, respectively.  
Then define the relation $\mathrel{\preccurlyeq}$ on ${P_{\mathrm w}}$ generated by 

(i)\ either $sa\mathrel{\preccurlyeq} sb$ when $ a\leq b$; or

(ii)\ $sabt\mathrel{\preccurlyeq} sct$ when $a\leq c$,

where $s,t\in{P_{\mathrm w}}$ and $a,b,c\in\mathbb N_{+}$. And let $u\mathrel{\preccurlyeq}\varepsilon$ for every $u\in{P_{\mathrm w}}$.
It is straightforward to verify that $\mathrel{\preccurlyeq}$ is a partial order on ${P_{\mathrm w}}$.

\begin{lemma}\label{lem:words-length-order}
Let $({P_{\mathrm w}},\mathrel{\preccurlyeq})$ be the poset defined above, and let $u,v\in{P_{\mathrm w}}$. If $u\mathrel{\preccurlyeq} v,$ then $ |u|\geq |v|.$
If $|u|=|v|=m\geq1$, then $u\mathrel{\preccurlyeq} v$ if and only if $u_i=v_i\ (1\leq i<m)$ and $ u_m\leq v_m.$
\end{lemma}
\begin{proof}
Let $u=w_0\mathrel{\preccurlyeq} w_1\mathrel{\preccurlyeq}\cdots\mathrel{\preccurlyeq} w_r=v$ be a chain of defining steps. By (i), (ii),
and $w\mathrel{\preccurlyeq}\varepsilon$,
$|u|=|w_0|\geq|w_1|\geq\cdots\geq|w_r|=|v|.$
If $|u|=|v|=m\geq1$, then
$|w_j|=m(0\leq j\leq r), w_j=p b_j, |p|=m-1, u_m=b_0\leq b_1\leq\cdots\leq b_r=v_m.$
Here $p=(u_1,\ldots,u_{m-1})=(v_1,\ldots,v_{m-1})$. This follows because only (i) preserves
positive length. Conversely, $u_i=v_i\ (1\leq i<m)$ and $u_m\leq v_m$
imply $u\mathrel{\preccurlyeq} v$ by (i).
\end{proof}

\begin{lemma}\label{lem:words-child-sup}
	Let $({P_{\mathrm w}},\mathrel{\preccurlyeq})$ be the poset defined above. For every $s\in{P_{\mathrm w}}$, put $C_s=\{sn:n\in\mathbb N_{+}\}$, then $C_s$ is a directed chain with $\bigvee C_s=s$.
\end{lemma}
\begin{proof}
	By (i), $s1\mathrel{\preccurlyeq} s2\mathrel{\preccurlyeq}\cdots$. Thus $C_s$ is a directed chain, because the larger of any two of its elements is an upper bound in $C_s$. Also $sn\mathrel{\preccurlyeq} s$ for every $n\in\mathbb N_{+}$:
	if $s=\varepsilon$, this follows from the maximality of $\varepsilon$; otherwise,
	writing $s=ta$, relation (ii) gives $tan\mathrel{\preccurlyeq} ta=s$.
	
	Let $\varnothing\ne E\subseteq C_s$ with $\{v(e):e\in E\}$ unbounded, and
	let $y$ be an upper bound of $E$. If $y=\varepsilon$, then $s\mathrel{\preccurlyeq} y$.
	Suppose $y\ne\varepsilon$ and choose $sn\in E$ with $n>v(y)$. By
Lemma~\ref{lem:words-length-order},
$sn\mathrel{\preccurlyeq} y\ \Longrightarrow\ |y|\leq |s|+1, |y|\ne |s|+1,$
so $|y|= |s|$. Thus $s$ is the least upper bound of $E$. Since
	$\{v(e):e\in C_s\}=\mathbb N_{+}$ is unbounded, $\bigvee C_s=s$.
\end{proof}

\begin{proposition}\label{prop:words-directed}
	Let $({P_{\mathrm w}},\mathrel{\preccurlyeq})$ be the poset defined above, every nonempty directed set $D\subseteq{P_{\mathrm w}}$ has a supremum.
That is, ${P_{\mathrm w}}$ is a countable dcpo.
\end{proposition}
\begin{proof}
Put $m=\min\{|d|:d\in D\}, D_m=\{d\in D:|d|=m\}.$ Choose $d_0\in D_m$. For any $d\in D$, there exists $e\in D$ with $d\mathrel{\preccurlyeq} e$ and $d_0\mathrel{\preccurlyeq} e$ by $D$ is a directed poset. By Lemma~\ref{lem:words-length-order}, $m\leq|e|\leq|d_0|=m$,
so $|e|=m$, that is, $e\in D_m$ and $d\mathrel{\preccurlyeq} e\in D_m$. Thus $D_m$ is a cofinal subset of the directed set $D$, and hence $D_m$ is directed.
	If $m=0$, then $\varepsilon\in D$, so $\varepsilon=\max D$.
Now suppose $m\geq 1$. By Lemma~\ref{lem:words-length-order}, any two elements of $D_m$
have the same first $m-1$ entries, so there exists a sequence $s$ with $|s|=m-1$ such that
$D_m\subseteq C_s$. If $\{v(d):d\in D_m\}$ is bounded, then
$D_m$ is a finite chain, hence has a greatest element $d^*$; since $D_m$ is a cofinal subset of
$D$, $d^*$ is also the greatest element of $D$, and $\bigvee D=d^*$. If $\{v(d):d\in D_m\}$ is unbounded, then	Lemma~\ref{lem:words-child-sup} applied to $E=D_m$ gives
$\bigvee D_m=s$; since $D_m$ is cofinal in $D$, directed suprema coincide, so $\bigvee D=\bigvee D_m=s$.
\end{proof}

\begin{lemma}\label{lem:words-grandchild}
Let $({P_{\mathrm w}},\mathrel{\preccurlyeq})$ be the partially ordered set defined above, every $ x\in{P_{\mathrm w}}$ and $i,j\in\mathbb N_{+}$, we have $xij\ll x$. 
\end{lemma}
\begin{proof}
We first consider the defining steps (i) and (ii). For $a\leq c$ and
$d\leq e$, both
$sabtd\mathrel{\preccurlyeq} sabte\mathrel{\preccurlyeq} scte,
 \ sabtd\mathrel{\preccurlyeq} sctd\mathrel{\preccurlyeq} scte$
are valid chains. The second places the step of type (ii) before the
step of type (i). Similarly, for $b\leq b'$ and $a\leq c$, the chain
$sab\mathrel{\preccurlyeq} sab'\mathrel{\preccurlyeq} sc$ can be replaced by $sab\mathrel{\preccurlyeq} sc$.
By successively applying these exchanges, a chain from $x$ to
$z=tb\ne\varepsilon$ may be chosen in the form
$
 x=w_0\mathrel{\preccurlyeq} w_1\mathrel{\preccurlyeq}\cdots\mathrel{\preccurlyeq} w_\ell=ta\mathrel{\preccurlyeq} tb=z,
 a\leq b,
$
where each step from $w_{r-1}$ to $w_r$ is of type (ii).
Appending $ij$ preserves steps of type (ii). Thus, for $k\geq j$,
$
 xij=w_0ij\mathrel{\preccurlyeq} w_1ij\mathrel{\preccurlyeq}\cdots\mathrel{\preccurlyeq} w_\ell ij
 =taij\mathrel{\preccurlyeq} tbj\mathrel{\preccurlyeq} tbk=zk.
$
In particular,
\begin{equation}\label{eq:words-grandchild-bound}
 xij\mathrel{\preccurlyeq} zk\qquad\text{whenever }x\mathrel{\preccurlyeq} z\ne\varepsilon
 \text{ and }k\geq j.
\end{equation}
For $z=\varepsilon$, repeated applications of (ii) with $c=a$,
followed by (i), give
\begin{equation}\label{eq:words-grandchild-empty}
 xij\mathrel{\preccurlyeq}\bigl(\operatorname{f}(xij)\bigr)\mathrel{\preccurlyeq}(k)
 \qquad\bigl(k\geq\operatorname{f}(xij)\bigr).
\end{equation}

Let $D\in\mathcal D^{\vee}({P_{\mathrm w}})$ and suppose $x\mathrel{\preccurlyeq}\bigvee D=z$.
If $D$ has a greatest element $d^*$, then
$
 xij\mathrel{\preccurlyeq} x\mathrel{\preccurlyeq} z=d^*\in D.
$
Otherwise, by Proposition~\ref{prop:words-directed}, there is a cofinal
set $E\subseteq D\cap C_z$ with $\{v(e):e\in E\}$ unbounded.
If $z\ne\varepsilon$, choose $zk\in E$ with $k\geq j$;
by \eqref{eq:words-grandchild-bound}, $xij\mathrel{\preccurlyeq} zk$.
If $z=\varepsilon$, choose $(k)\in E$ with
$k\geq\operatorname{f}(xij)$;
by \eqref{eq:words-grandchild-empty}, $xij\mathrel{\preccurlyeq}(k)$.
In each case $D\cap\mathord{\uparrow}xij\ne\varnothing$.
Since $D$ was arbitrary, $xij\ll x$.
\end{proof}

\begin{theorem}\label{thm:words-nondirected}
	Let $({P_{\mathrm w}},\mathrel{\preccurlyeq})$ be the poset defined above, the poset ${P_{\mathrm w}}$ is a countable strongly meet continuous dcpo, and ${P_{\mathrm w}}$ is not continuous.
\end{theorem}
\begin{proof}
By Proposition~\ref{prop:words-directed}, ${P_{\mathrm w}}$ is a countable dcpo.
Let $x\in{P_{\mathrm w}}$ and $H\in\Gamma({P_{\mathrm w}})$ with $\mathord{\Downarrow} x\subseteq H$.
By Lemma~\ref{lem:words-grandchild},
$\{xij:i,j\in\mathbb N_{+}\}\subseteq\mathord{\Downarrow} x\subseteq H.$
For each $i\in\mathbb N_{+}$, Lemma~\ref{lem:words-child-sup} and $H\in\Gamma({P_{\mathrm w}}),$ then
$ xi=\bigvee_{j\in\mathbb N_{+}}xij\in H.$
And we have
$x=\bigvee_{i\in\mathbb N_{+}}xi\in H.$
As $H$ was arbitrary,
$x\in\bigcap\{H\in\Gamma({P_{\mathrm w}}):\mathord{\Downarrow} x\subseteq H\}
 =\operatorname{cl}_{\sigma}(\mathord{\Downarrow} x).$
Hence ${P_{\mathrm w}}$ is strongly meet continuous by
Theorem~\ref{thm:smc-equivalences}.

Fix $x\in{P_{\mathrm w}}\setminus\{\varepsilon\}$ and put $m=|x|\geq1$.
For $u\ll x$, Lemma~\ref{lem:words-child-sup} gives
$
 x=\bigvee_{n\in\mathbb N_{+}}xn,
 \ u\mathrel{\preccurlyeq} xn\ \text{for some }n\in\mathbb N_{+}.
$
It follows from Lemma~\ref{lem:words-length-order} that
$|u|\geq m+1$.
Suppose $|u|=m+1$. By Lemma~\ref{lem:words-length-order}, $u=xr$ for some $r\in\mathbb N_{+}$.
$
 x^+=(x_1,\ldots,x_{m-1},v(x)+1).
$
Since $x\mathrel{\preccurlyeq} x^+=\bigvee_{n\in\mathbb N_{+}}x^+n$ and $u\ll x$,
there exists $n\in\mathbb N_{+}$ such that
$
 xr=u\mathrel{\preccurlyeq} x^+n,
  |xr|=|x^+n|=m+1.
$
Lemma~\ref{lem:words-length-order} now yields $x=x^+$,
a contradiction. We have proved
\begin{equation}\tag{$\dagger$}
 |u|\geq m+2\ \text{for every }u\ll x.
\end{equation}

Put $a=x11$ and $b=x21$. By Lemma~\ref{lem:words-grandchild},
$a,b\in\mathord{\Downarrow} x$. Suppose that
$c\in\mathord{\Downarrow} x\cap\mathord{\uparrow} a\cap\mathord{\uparrow} b$. By $(\dagger)$ and
Lemma~\ref{lem:words-length-order},
$m+2\leq|c|\leq\min\{|a|,|b|\}=m+2.$
Thus $|a|=|b|=|c|=m+2$. And $a=x11\mathrel{\preccurlyeq} c, b=x21\mathrel{\preccurlyeq} c.$
But by Lemma~\ref{lem:words-length-order},
$(c_1,\ldots,c_{m+1})=x1,
 \ (c_1,\ldots,c_{m+1})=x2,$
a contradiction. Therefore
$
 \mathord{\Downarrow} x\cap\mathord{\uparrow} a\cap\mathord{\uparrow} b=\varnothing,
 \forall\ a,b\in\mathord{\Downarrow} x.
$
Consequently, $\mathord{\Downarrow} x$ is not directed and ${P_{\mathrm w}}$ is not continuous.
\end{proof}

\begin{proposition}\label{prop:smc-properties}
Let $P$ be a strongly meet continuous poset, then for every $x\in P$, $\bigvee\mathord{\Downarrow} x=x$ and $\mathord{\Downarrow} x$ is Scott irreducible.
\end{proposition}
\begin{proof}
Fix $x\in P$. Since $P$ is strongly meet continuous, $x\in\operatorname{cl}_{\sigma}(\mathord{\Downarrow} x)$. If $y$ is an upper bound of $\mathord{\Downarrow} x$,
then $\mathord{\Downarrow} x\subseteq\mathord{\downarrow} y$. Since $\mathord{\downarrow} y$ is Scott closed,
$\mathord{\downarrow} x=\operatorname{cl}_{\sigma}(\mathord{\Downarrow} x)\subseteq\mathord{\downarrow} y,$
and therefore $x\leq y$. Then $\bigvee\mathord{\Downarrow} x=x$.

Let $U,V\in\sigma(P)$ both meet
$\mathord{\Downarrow} x$, and choose $a\in U\cap\mathord{\Downarrow} x$ and $b\in V\cap\mathord{\Downarrow} x$. Since $U$ and $V$ are
upper sets and $a,b\leq x$, we have $x\in U\cap V$. Then $U\cap V\in\sigma(P)$, so
Theorem~\ref{thm:smc-equivalences}(iv) yields
$U\cap V\cap\mathord{\Downarrow} x\ne\varnothing$. Hence $\mathord{\Downarrow} x$ is irreducible.
\end{proof}

\begin{theorem}\label{prop:open-union}
Let $P$ be a poset. $P$ is strongly meet continuous if and only if
 $\forall\ U\in\sigma(P), U=\bigcup_{a\in U}\mathord{\Uparrow} a.$
\end{theorem}
\begin{proof}
The inclusion $\bigcup_{a\in U}\mathord{\Uparrow} a\subseteq U$ is immediate from the definition of
$\mathord{\Uparrow} a$. The reverse inclusion states exactly that every $x\in U$ satisfies
$U\cap\mathord{\Downarrow} x\ne\varnothing$, which is condition (iv) of
Theorem~\ref{thm:smc-equivalences}. 
\end{proof}

\section{The relation between continuous and strongly meet continuous}\label{sec:continuous-case}

We now investigate when a strongly meet continuous poset is continuous and establish
several equivalent characterizations. 

Let $P$ be any poset and $A\subseteq P$. Define 
$$A':=\{x\in P:D\in\mathcal D^{\vee}(P),\ D\subseteq\mathord{\downarrow} A. x=\vee D\},$$
$$A'':=\{x\in P:\exists D\in\mathcal D^{\vee}(P),\ D\subseteq\mathord{\downarrow} A,\ x\leq\bigvee D\}.$$
Then a poset $P$ has \emph{one step closure}\cite{ref4} if
$A'=\operatorname{cl}_{\sigma} A$ for every $A\subseteq P$, and \emph{weak one step closure}\cite{ref5} if $A''=\operatorname{cl}_{\sigma} A$ for every $A\subseteq P$. 

\begin{proposition}[\textup{\cite{ref5}}]\label{prop:4.1}
A poset $P$ has one step closure if and only if it is meet continuous and has weak one
step closure.
\end{proposition}

\begin{theorem}[\textup{\cite{ref4}}]\label{thm:4.2}
Every continuous poset has one step closure.
\end{theorem}

\begin{example}[\textup{\cite{ref5}}]\label{ex:sorgenfrey}
Let $\mathbb R_\ell$ denote the Sorgenfrey line, and let $K(\mathbb R_\ell)$ be the poset
of its nonempty compact saturated subsets ordered by reverse inclusion. Then
$K(\mathbb R_\ell)$ has one step closure but is not continuous.
\end{example}

The main result of this subsection shows that, under strongly meet continuity, either one
step closure or weak one step closure is equivalent to continuity.

\begin{theorem}\label{thm:4.4}
For a poset $P$, the following conditions are equivalent:
\begin{enumerate}[label=(\roman*)]
 \item $P$ is continuous;
 \item $P$ is strongly meet continuous and has one step closure;
 \item $P$ is strongly meet continuous and has weak one step closure;
 \item $x$ is the supremum of a directed subset of $\mathord{\Downarrow} x$ for every $x\in P$;
 \item $x\in(\mathord{\Downarrow} x)'$ for every $x\in P$.
\end{enumerate}
\end{theorem}
\begin{proof}
(i)$\Rightarrow$(ii). This follows from Proposition~\ref{prop:continuous-smc-mc} and
Theorem~\ref{thm:4.2}.

(ii)$\Rightarrow$(iii). If $P$ has one step closure, then it has weak one step closure,
since $A''=\mathord{\downarrow}A'=\mathord{\downarrow}\operatorname{cl}_{\sigma} A=\operatorname{cl}_{\sigma} A$ for every $A\subseteq P$.

(iii)$\Rightarrow$(iv). Fix $x\in P$. By
Theorem~\ref{thm:smc-equivalences}(iii) we have $\operatorname{cl}_{\sigma}(\mathord{\Downarrow} x)=\mathord{\downarrow} x$, so
$x\in(\mathord{\Downarrow} x)''$. By the definition of the operator $''$, there exists a directed set
$E_x\subseteq\mathord{\downarrow}\mathord{\Downarrow} x=\mathord{\Downarrow} x$ such that $x\leq\bigvee E_x$. On the other hand,
$E_x\subseteq\mathord{\Downarrow} x\subseteq\mathord{\downarrow} x$ implies $\bigvee E_x\leq x$. Hence $\bigvee E_x=x$,
which proves (iv).

(iv)$\Rightarrow$(i). This is exactly Lemma~\ref{lem:directed-approx}.

(iv)$\Leftrightarrow$(v). Since $\mathord{\Downarrow} x$ is a lower set, we have
$\mathord{\downarrow}\mathord{\Downarrow} x=\mathord{\Downarrow} x$ and the elements of $(\mathord{\Downarrow} x)'$ are precisely the suprema of
directed subsets of $\mathord{\Downarrow} x$.
\end{proof}

We can see the conjunction of strongly meet continuity and one step closure characterizes
continuity. The implications established so far are summarized in the following diagram.

\begin{center}
\begin{tikzpicture}[
 >=Stealth,
 concept/.style={draw,rounded corners,align=center,inner sep=5pt,
 minimum width=2.9cm,minimum height=9mm},
 merge/.style={circle,fill=black,inner sep=1.3pt}
 ]
 \node[concept] (cont) at (-3.6,0) {Continuous};
 \node[concept] (smc) at (0,1.6) {Strongly\\meet continuous};
 \node[concept] (osc) at (0,-1.6) {One step\\closure};
 \node[concept] (mc) at (3.6,0) {Meet continuous};
 \draw[->] (cont) -- (smc);
 \draw[->] (cont) -- (osc);
 \draw[->] (smc) -- (mc);
 \draw[->] (osc) -- (mc);
 \draw (smc.south) -- (osc.north);
 \node[merge] (join) at (0,0) {};
 \draw[->] (join.west) -- (cont.east);
\end{tikzpicture}
\end{center}

\begin{proposition}[\textup{\cite{ref4}}]\label{prop:4.1prop}
Let $P$ be a poset with one step closure. Then $\mathord{\Uparrow} a=\operatorname{int}_{\sigma}(\mathord{\uparrow} a)$ for each $ a\in P$.
\end{proposition}

By Proposition~\ref{prop:4.1prop}, under one step closure the way above set $\mathord{\Uparrow} a$
coincides with the Scott interior of the principal filter $\mathord{\uparrow} a$. This suggests examining
the Scott openness of $\mathord{\Uparrow} a$ and the interpolation property of $\ll$ under strongly meet
continuity.

\begin{theorem}\label{prop:4.5}
Let $P$ be a strongly meet continuous poset, the following conditions are equivalent:
\begin{enumerate}[label=(\roman*)]
 \item $P$ is continuous;
 \item $\mathord{\Downarrow} x$ is directed for every $x\in P$;
 \item $\mathord{\Uparrow} a\in\sigma(P)$ for every $a\in P$;
 \item $\ll$ has the interpolation property;
 \item $\displaystyle\mathord{\Downarrow} x=\bigcup_{b\ll x}\mathord{\Downarrow} b$ for every $x\in P$;
 \item $x\in\operatorname{int}_{\sigma}(\mathord{\uparrow} a\cap\mathord{\uparrow} b)$ if $a,b\ll x$.
\end{enumerate}
\end{theorem}
\begin{proof}
(i)$\Rightarrow$(ii) holds by the definition of continuity, and (i)$\Rightarrow$(iv)
follows from Lemma~\ref{lem:interpolation}.

(iv)$\Rightarrow$(iii). Fix $a\in P$. The set $\mathord{\Uparrow} a$ is an upper set by the
transitivity of $\ll$. Let $D\in\mathcal D^{\vee}(P)$ with $\bigvee D\in\mathord{\Uparrow} a$, that is,
$a\ll\bigvee D$. By interpolation there exists $b\in P$ with $a\ll b\ll\bigvee D$, and
then there exists $d\in D$ with $b\leq d$. It follows that $a\ll d$, so
$d\in\mathord{\Uparrow} a$ and $D\cap\mathord{\Uparrow} a\ne\varnothing$. Hence $\mathord{\Uparrow} a\in\sigma(P)$.

(iii)$\Rightarrow$(ii). Let $x\in P$ and $a,b\in\mathord{\Downarrow} x$. Then $\mathord{\Uparrow} a$ and $\mathord{\Uparrow} b$ are
Scott open neighborhoods of $x$, so $\mathord{\Uparrow} a\cap\mathord{\Uparrow} b\in\sigma(P)$ contains $x$. By
Theorem~\ref{thm:smc-equivalences}(iv) there exists $c\in\mathord{\Downarrow} x$ with
$c\in\mathord{\Uparrow} a\cap\mathord{\Uparrow} b$, that is, $a\ll c$ and $b\ll c$; then $a,b\leq c$. Hence
$\mathord{\Downarrow} x$ is directed.

(ii)$\Rightarrow$(i). This follows from Proposition~\ref{prop:smc-properties}, which
states that $\bigvee\mathord{\Downarrow} x=x$ for every $x\in P$: together with (ii) this is exactly the
definition of continuity.

(iv)$\Leftrightarrow$(v). The inclusion $\bigcup_{b\ll x}\mathord{\Downarrow} b\subseteq\mathord{\Downarrow} x$ always holds, while the reverse inclusion says precisely that $a\ll x$ implies $a\ll b\ll x$ for some $b$.

(iii)$\Rightarrow$(vi). If $a,b\ll x$, then $x\in\mathord{\Uparrow} a\cap\mathord{\Uparrow} b$ and $\mathord{\Uparrow}
a\cap\mathord{\Uparrow} b$ is Scott open. Since $\mathord{\Uparrow} a\cap\mathord{\Uparrow} b\subseteq\mathord{\uparrow} a\cap\mathord{\uparrow} b$, we obtain
$x\in\operatorname{int}_{\sigma}(\mathord{\uparrow} a\cap\mathord{\uparrow} b)$.

(vi)$\Rightarrow$(ii). Let $a,b\in\mathord{\Downarrow} x$. By (vi),
$\operatorname{int}_{\sigma}(\mathord{\uparrow} a\cap\mathord{\uparrow} b)$ is a Scott open neighborhood of $x$, so
Theorem~\ref{thm:smc-equivalences}(iv) provides $c\in\mathord{\Downarrow} x\cap\mathord{\uparrow} a\cap\mathord{\uparrow} b$.
Then $a\leq c$, $b\leq c$ and $c\ll x$, which proves that $\mathord{\Downarrow} x$ is directed.
\end{proof}

\begin{lemma}\label{lem:4.7}
If a nonempty lower directed set $F$ is covered by finitely many upper sets
$U_1,\ldots,U_m$, then $F\subseteq U_j$ for some $j$.
\end{lemma}
\begin{proof}
If every $1\leq j\leq m, F\not\subseteq U_j$. For each $1\leq j\leq m$, choose $z_j\in F\setminus U_j$, and then
choose $z\in F$ with $z\leq z_j$ for every $j$; such a $z$ exists because $F$ is lower
directed. Since the $U_j$ cover $F$, we have $z\in U_j$ for some $j$. As $U_j$ is an
upper set and $z\leq z_j$, this implies $z_j\in U_j$, contradicting the choice of $z_j$.
\end{proof}

\begin{lemma}\label{lem:4.8}
Let $P$ be a poset. If $\forall\ x\in U\in\sigma(P)$, there exists $a\in U$ such that
$x\in\operatorname{int}_{\sigma}(\mathord{\uparrow} a)$, then $P$ is continuous.
\end{lemma}
\begin{proof}
For each $x\in P$, put $B_x=\{a\in P:x\in\operatorname{int}_{\sigma}(\mathord{\uparrow} a)\}.$
Since $\operatorname{int}_{\sigma}(\mathord{\uparrow} a)\in\sigma(P)$ and $\operatorname{int}_{\sigma}(\mathord{\uparrow} a)\subseteq\mathord{\uparrow} a$, every
$a\in B_x$ satisfies $x\in\mathord{\uparrow} a$; moreover, if $D\in\mathcal D^{\vee}(P)$ and
$x\leq\bigvee D$, then there exists $d\in D\cap \operatorname{int}_{\sigma}(\mathord{\uparrow} a)$, and
$a\leq d$. This shows that $B_x\subseteq\mathord{\Downarrow} x$. Taking $U=P$ in the hypothesis, we obtain
$B_x\ne\varnothing$.

Let $a,b\in B_x$. Then
$x\in\operatorname{int}_{\sigma}(\mathord{\uparrow} a)\cap\operatorname{int}_{\sigma}(\mathord{\uparrow} b)\in\sigma(P)$. By the hypothesis, there exists $c\in\operatorname{int}_{\sigma}(\mathord{\uparrow}
a)\cap\operatorname{int}_{\sigma}(\mathord{\uparrow} b)$ with $x\in\operatorname{int}_{\sigma}(\mathord{\uparrow} c)$. Then $c\in B_x$ and $c\in\mathord{\uparrow} a\cap\mathord{\uparrow}
b$, that is, $a\leq c$ and $b\leq c$. Hence $B_x$ is directed.

Finally, we show that $\bigvee B_x=x$. Clearly $x$ is an upper bound of $B_x$. If $y$ is
an upper bound of $B_x$ with $x\not\leq y$, then $x\in P\setminus\mathord{\downarrow} y\in\sigma(P)$,
then there exists $a\in B_x$ with $a\not\leq y$, contradicting the assumption that $y$ is an upper bound of $B_x$. Therefore $\bigvee B_x=x$. Since $B_x\subseteq\mathord{\Downarrow} x$ is directed
with supremum $x$, Lemma~\ref{lem:directed-approx} shows that $P$ is continuous.
\end{proof}

\begin{lemma}\label{lem:scott-open-filter-base}
	If $P$ is strongly meet continuous, then $\sigma(P)$ has a basis consisting of Scott open
	filters.
\end{lemma}
\begin{proof}
	Let $x\in U\in\sigma(P)$. We construct inductively a sequence $\{a_n:n\geq 0\}$  by
	Theorem~\ref{thm:smc-equivalences}(iv) such
	that $a_0=x$ and $a_{n+1}\in U\cap\mathord{\Downarrow}a_n$ for every $n\geq 0$. Put
	$F=\bigcup_{n\geq 0}\mathord{\uparrow}a_n.$
	We verify that $F$ is a Scott open filter with $x\in F\subseteq U$.
	Since $x=a_0\in\mathord{\uparrow}a_0\subseteq F$ and each $a_n\in U$ with $U$ an upper set, we have
	$x\in F\subseteq U$. Being a union of upper sets, $F$ is an upper set. To see that $F$ is
	lower directed, let $u,v\in F$ and choose $m,n\geq 0$ with $a_m\leq u$ and $a_n\leq v$.
	Since $a_{k+1}\ll a_k$ implies $a_{k+1}\leq a_k$, the sequence $\{a_n\}$ is decreasing.
	Taking $k=\max\{m,n\}$, we obtain $a_k\leq a_m\leq u$ and $a_k\leq a_n\leq v$, so $a_k$ is
	a common lower bound of $u$ and $v$ lying in $F$. 
	Finally, let $D\in\mathcal D^{\vee}(P)$ with $\bigvee D\in F$. Then $a_n\leq\bigvee D$ for some
	$n\geq 0$. Since $a_{n+1}\ll a_n\leq\bigvee D$, there exists $d\in D$ with
	$a_{n+1}\leq d$, and hence $d\in\mathord{\uparrow}a_{n+1}\subseteq F$. Thus $D\cap F\ne\varnothing$, and
	$F$ is Scott open.
\end{proof}

\begin{theorem}\label{prop:4.9}
Let $P$ be a poset, then $P$ is continuous if and only if $\Sigma P$ is core compact and
$\sigma(P)$ has a basis consisting of Scott open filters.
\end{theorem}
\begin{proof}
Suppose that $P$ is continuous. Let $x\in U\in\sigma(P)$ and choose
$a\in U\cap\mathord{\Downarrow} x$ by Lemma~\ref{lem:interpolation}. Then
$x\in\mathord{\Uparrow} a\subseteq\mathord{\uparrow} a\subseteq U.$
Moreover, $\mathord{\Uparrow} a\ll U$: indeed, if $\mathcal W$ is an open cover of
$U$, there exists $W\in\mathcal W$ containing $a$, and since $W$ is an upper set,
$\mathord{\uparrow} a\subseteq W$; it follows that $\mathord{\Uparrow} a\subseteq W$. Hence $\Sigma P$ is core
compact. By Proposition~\ref{prop:continuous-smc-mc} and
Lemma~\ref{lem:scott-open-filter-base}, $\sigma(P)$ has a basis consisting of Scott open
filters.
Conversely, suppose that $\Sigma P$ is core compact and that $\sigma(P)$ has a basis
consisting of Scott open filters. Let $x\in U\in\sigma(P)$. By core compactness, choose
$V\in\sigma(P)$ with $x\in V\ll U$, and then choose a Scott open
filter $F$ with $x\in F\subseteq V$. We claim that there exists $a\in U$ with $F\subseteq\mathord{\uparrow} a$. Suppose not. Then for each
$a\in U$, there exists $z_a\in F$ with $a\not\leq z_a$, and the sets $O_a=P\setminus\mathord{\downarrow}z_a, a\in U,$
are Scott open. For each $b\in U$ we have
$b\not\leq z_b$, so $b\in O_b$; thus $\{O_a:a\in U\}$ is an open cover of $U$. Since
$V\ll U$ and $F\subseteq V\subseteq U$, there exist
$a_1,\ldots,a_m\in U$ such that $F\subseteq\bigcup_{j=1}^m O_{a_j}$. By
Lemma~\ref{lem:4.7} applied to the lower directed set $F$ and the upper sets
$O_{a_1},\ldots,O_{a_m}$, we have $F\subseteq O_{a_j}$ for some $j$; but
$z_{a_j}\in F\setminus O_{a_j}$, a contradiction. 
Thus, there exists $a\in U$ such that $x\in F\subseteq\mathord{\uparrow} a\subseteq U.$
Since $F$ is open and $F\subseteq\mathord{\uparrow} a$, we have $x\in F\subseteq\operatorname{int}_{\sigma}(\mathord{\uparrow} a)$. Thus  $P$ is continuous by Lemma~\ref{lem:4.8}.
\end{proof}

\begin{corollary}\label{cor:4.10}
A poset $P$ is continuous if and only if $P$ is strongly meet continuous and $\Sigma P$ is
core compact.
\end{corollary}
\begin{proof}
If $P$ is continuous, then $P$ is strongly meet continuous by
Proposition~\ref{prop:continuous-smc-mc}, and $\Sigma P$ is core compact as shown in the
proof of Theorem~\ref{prop:4.9}. The converse follows from
Lemma~\ref{lem:scott-open-filter-base} and Theorem~\ref{prop:4.9}.
\end{proof}

\begin{corollary}\label{cor:smc-corecompact-sober}
Let $P$ be a strongly meet continuous dcpo. If $\Sigma P$ is core compact, then $\Sigma P$ is sober.
\end{corollary}
\begin{proof}
By Corollary~\ref{cor:4.10}, $P$ is continuous. Since the Scott space of every continuous dcpo is sober \cite{ref11}, $\Sigma P$ is sober.
\end{proof}

\begin{proposition}[\textup{\cite{ref21}}]\label{prop:principal-ideal-continuity}
Let $P$ be a dcpo, then $P$ is continuous if and only if every principal ideal $\mathord{\downarrow} z$,
$z\in P$, is continuous.
\end{proposition}

\begin{corollary}\label{cor:4.11}
Let $P$ be a strongly meet continuous dcpo, then $P$ is continuous if and only if
$\Sigma(\mathord{\downarrow} z)$ is core compact for every $z\in P$. Moreover, if every two elements
of each principal ideal $\mathord{\downarrow} z$ have a join in that ideal, then $P$ is continuous.
\end{corollary}
\begin{proof}
Fix $z\in P$ and put $C=\mathord{\downarrow} z$. For $U\in\sigma(C)$, define
$\widetilde U=U\cup(P\setminus C)$. Then $\widetilde U\in\sigma(P)$: it is clearly an
upper set, and if $D\in\mathcal D^{\vee}(P)$ with $\bigvee D\in\widetilde U$, then either
$\bigvee D\in P\setminus C$, in which case $D$ meets the Scott open set $P\setminus C$, or
$\bigvee D\in U\subseteq C$, in which case $D\subseteq\mathord{\downarrow}\bigvee D\subseteq C$ and
$D$ meets $U$ because $U\in\sigma(C)$. Let $x\in U\in\sigma(C)$. Since $x\in\widetilde U\in\sigma(P)$ and $P$ is strongly
meet continuous, Theorem~\ref{thm:smc-equivalences}(iv) gives $a\in\widetilde U$ with
$a\ll_P x$. As $a\leq x\leq z$, we have $a\in C$, hence $a\in U$; and every directed
subset of $C$ is directed in $P$ with the same supremum, so $a\ll_C x$. By
Theorem~\ref{thm:smc-equivalences}, $C$ is strongly meet continuous. Applying
Proposition~\ref{prop:principal-ideal-continuity} and Corollary~\ref{cor:4.10}, we
conclude that
$P\text{ is continuous}$ if and only if $\forall z\in P,\ \mathord{\downarrow} z\text{ is continuous}$ if and only if $ \forall z\in P,\ \Sigma(\mathord{\downarrow} z)\text{ is core compact}.$

Suppose that every two elements of each principal ideal $C=\mathord{\downarrow}
z$ have a join in $C$. Let $x\in C$. The set $\mathord{\Downarrow}_C x$ is nonempty. Given
$a,b\in\mathord{\Downarrow}_C x$, their join $a\vee_C b$ exists in $C$, and $a\vee_C b\in\mathord{\Downarrow}_C x$: if $D\subseteq C$ is directed with $x\leq\bigvee D$, then there
exist $d_1,d_2\in D$ with $a\leq d_1$ and $b\leq d_2$, and directedness gives $d\in D$
with $d_1,d_2\leq d$, so $a\vee_C b\leq d$. Hence $\mathord{\Downarrow}_C x$ is directed, and since
$\bigvee\mathord{\Downarrow}_C x=x$, the ideal $C$ is continuous. By
Proposition~\ref{prop:principal-ideal-continuity}, $P$ is continuous.
\end{proof}

\section{Conclusion}\label{sec:conclusion}

In this paper, we have introduced strongly meet continuous posets and characterized them by
$x\in\operatorname{cl}_{\sigma}(\Downarrow\!x)$.
The countable dcpo constructed in this paper shows that strongly meet continuity
does not imply continuity. For strongly
meet continuous posets, continuity is equivalent to either one step closure or
weak one step closure. Every strongly meet continuous poset also admits a basis of Scott open
filters; together with core compactness, this gives a topological
characterization of continuous posets. Thus strongly meet continuity provides a
link between order theoretic approximation, one step closure and the topology
of Scott spaces. In future work, it would be interesting to identify conditions under which a strongly meet continuous dcpo has a sober Scott space without assuming core compactness. 

\section*{Acknowledgments}
This work was supported by the National Natural Science Foundation of China
(Grant No.~12331016).

\section*{Declaration of AI}
The authors declare that no AI tools were used in this manuscript.

\end{document}